\RequirePackage{fix-cm}
\documentclass[smallextended]{svjour3}       
\smartqed  
\usepackage{amsmath}%
\usepackage{amsfonts}%
\usepackage{graphicx}
\usepackage{natbib}
\begin{document}

\title{A Property of the Kullback--Leibler Divergence for Location-scale Models
}


\author{Cristiano Villa
}


\institute{C. Villa \at
              School of Mathematics, Statistics and Actuarial Science, University of Kent, Canterbury, UK \\
              \email{cv88@kent.ac.uk}           
}

\date{Received: date / Accepted: date}

\maketitle

\begin{abstract}
In this paper, we discuss a property of the Kullback--Leibler divergence measured between two models of the family of the location-scale distributions. We show that, if model $M_1$ and model $M_2$ are represented by location-scale distributions, then the minimum Kullback--Leibler divergence from $M_1$ to $M_2$, with respect to the parameters of $M_2$, is independent from the value of the parameters of $M_1$. Furthermore, we show that the property holds for models that can be transformed into location-scale distributions. We illustrate a possible application of the property in objective Bayesian model selection.
\keywords{Bayes factor \and Kullback--Leibler divergence \and Location-scale model \and Model prior \and Objective Bayes \and Self-information loss}
 \subclass{62C05 \and 62B10 \and 62F03}
\end{abstract}

\section{Introduction and notation}\label{sc_intro}
There are various circumstances in statistics where the Kullback--Leibler divergence \citep{Kull:1951} can be applied. For example, the well known Akaike Information Criterion \citep{Aka:1973} for model selection is based on the relationship between the Kullback--Leibler divergence and the maximized likelihood; relationship that will be exploited in the paper. Other applications of the Kullback--Leibler divergence in model selection problems have been discussed by \cite{Zheng:2004} and \cite{LvLiu:2014}. In the former paper, the authors perform the model selection on the basis of the ratio between two Kullback--Leibler divergences, where each competing model is compared with the true (and unknown) model. In the latter paper, a generalised version of the Bayesian information criterion for Generalized Linear Models is introduced by taking into consideration model misspecification.

More in general. Statistical inference starts with a set of observations $x=(x_1\ldots,x_n)$, which is assumed to have been generated by an unknown model, say $f(x)$. By means of estimation methods, either frequentist or Bayesian, an approximated model, say $g(x)$, is obtained. The overall aim is to have $g(x)$ as ``similar'' as possible to $f(x)$. The Kullback--Leibler divergence between the \emph{true} model $f(x)$ and the \emph{approximated} model $g(x)$ is given by
\begin{eqnarray}\label{eq_intro_1}
D_{KL}(f(x)\|g(x)) &=& \int_{\mathcal{X}} f(x)\log\{f(x)/g(x)\}\;dx \nonumber \\
&=& \mathbb{E}_f[\log f(x)] - \mathbb{E}_f[\log g(x)],
\end{eqnarray}
where $\mathcal{X}$ is the support of the true model, and the expectations are taken with respect to $f(x)$. The divergence in \eqref{eq_intro_1} measures the dissimilarity between two distributions, and it can be interpreted as the loss in information when $g(x)$ is used to approximate the true model $f(x)$. Therefore, the idea is to chose $g(x)$ which minimizes \eqref{eq_intro_1}. Given that the true model is unknown, in practical situations only the second expectation on the right-had-side of \eqref{eq_intro_1} has to be considered for model comparison or inference. In addition, the approximated model would be defined up to an unknown parameter (or vector of parameters) $\theta$, which carries additional uncertainty. In a frequentist set up, $\theta$ can be estimated through the maximum likelihood principle, whilst in the Bayesian framework one defines a prior distribution, say $\pi(\theta)$, which represents the uncertainty about the true parameter value, and obtains a posterior distribution of the parameter given the observations.
Should we be interested in performing model selection between two distributions (or testing two alternative hypothesis), both models will be defined up to some parameters, $f(x;\eta)$ and $g(x;\theta)$. In this case, there is a further level of uncertainty carried by $\eta$ itself, and the Kullback--Leibler divergence will be of the form
\begin{eqnarray}\label{eq_intro_2}
D_{KL}(f(x;\eta)\|g(x;\theta)) &=& \int_{\mathcal{X}} f(x;\eta)\log\{f(x;\eta)/g(x;\theta)\}\;dx \nonumber \\
&=& \mathbb{E}_f[\log f(x;\eta)] - \mathbb{E}_f[\log g(x;\theta)],
\end{eqnarray}
In this case, the uncertainty carried by $\eta$ implies that the inferential goal is to find the most similar model to $f(x;\eta)$ in expectation. In the Bayesian framework the expectation would be with respect to a prior assigned to $\eta$, representing the uncertainty about its true value.

In this paper we show that, when $f(x;\eta)$ and $g(x;\theta)$ belong to the location-scale family of distributions, the minimum Kullback--Leibler in \eqref{eq_intro_2} with respect to $\theta$, does not depend on $\eta$. As such, the loss in information is choosing $g(\cdot)$ as an approximation to $f(\cdot)$, is not affected by the uncertainty carried by $\eta$. The result is not limited to distributions of the location-scale family, but it is extended to distributions that, through a change in the variable, can be transformed into location-scale distributions. We illustrate the above property by an application to objective Bayesian model priors.

The paper is organized as follows. In Section \ref{sc_locsca} we introduce the family of location-scale distributions. Section \ref{sc_minkl} shows the main result, where we show how the minimum Kullback--Leibler divergence between any two location-scale models does not depend on the parameter of the first model, that is, the model from which the divergence is measured. We illustrate the result with two examples. An example of the usefulness of the main result in Bayesian model selection is shown in Section \ref{sc_modsel}. Finally, Section \ref{sc_disc} is dedicated to final remarks and discussions.

\section{Location-scale models}\label{sc_locsca}
We start by introducing the family of location-scale distributions. A detailed discussion of this family can be found, for example, in \cite{John:1992} and \cite{John:1994}.

We say that random variable $X$ belongs to the location-scale family if its distribution is of the form
$$F_X(x;a,b)=F\left(\frac{x-a}{b}\right) \qquad x\in\mathcal{X}\subseteq\mathbb{R}, a\in\mathbb{R},b>0,$$
where $F(\cdot)$ is a distribution function which characterizes the member of the family; $a$ and $b$ are, respectively, the location parameter and the scale parameter of the distribution. The location parameter $a$ can take any real value and it represents the position of the model on the abscissa, whilst the scale parameter $b$ can take positive values and represents the dispersion of the model. Location-scale distributions do not have any other parameter: that is, they are completely defined by $a$ and $b$.

Assuming that the distribution function of the random variable $X$ is absolutely continuous, the corresponding density function is
$$f(x;a,b) = b^{-1}h\left(\frac{x-a}{b}\right),$$
where $h(\cdot)$ is called the reduced density of $X$ and it has location zero and scale one.

If we consider two location-scale distributions with different reduced density, but same location and scale parameters, then they differ in the shape only - i.e. skewness and kurtosis. On the other hand, two location-scale distributions with the same reduced density are said to belong to the same type of location-scale distributions, which differ only in the values of the parameters.

The distribution functions belonging to the location-scale family can be categorized in two groups: genuine location-scale distributions, and distributions that require transformations to become location-scale models. The former group includes, for example, the normal distribution, the two-parameters exponential distribution, the extreme value distributions of type I for the maximum and the minimum, the logistic distribution and the uniform distribution. Among the distributions that can be transformed into location-scale models we have, for example, the translated Weibull distribution, where the transformation $\tilde{x}=\log(x-a)$ generates a location-scale type I extreme value distribution for the minimum with parameters $\tilde{b}=1/c$ and $\tilde{a}=\log b$. The log-normal distribution, where the transformation $\tilde{x}=\log(x-a)$ generates a normal with mean $\tilde{a}=a$ and standard deviation $\tilde{b}=b$.

In addition, the family includes distributions where either the location or the scale parameter can be omitted. In particular, by setting $b=1$, we have a location distribution, and by setting $a=0$, we have a scale distribution. Well known examples of the second type are the exponential and the half-normal densities.

\section{Minimum Kullback--Leibler divergence between two location-scale models}\label{sc_minkl}
Let us consider the densities $f_1(x;\eta)$ and $f_2(x;\theta)$, where both $\eta$ and $\theta$ are unknown. The Kullback--Leibler divergence between the two distributions is given by
\begin{eqnarray}\label{eq_kl_1}
D_{KL}(f_1(x;\eta)\|f_2(x;\theta)) &=& \int_{\mathcal{X}} f_1(x;\eta)\log\{f_1(x;\eta)/f_2(x;\theta)\}\;dx \nonumber \\
&=& \mathbb{E}_{f_1}[\log f_1(x;\eta)] - \mathbb{E}_{f_1}[\log f_2(x;\theta)].
\end{eqnarray}
Let us assume that we have observed data generated by the unknown model $f_1(x;\theta)$, say $x=(x_1,\ldots,x_n)$. Then, the sample Kullback--Leibler divergence between the densities is given by
$$D_{KL}^n = n^{-1}\sum_{i=1}^n\log\{f_1(x_i;\eta)/f_2(x_2;\theta)\},$$
which, as $n\rightarrow\infty$, converges to \eqref{eq_kl_1} with probability one. This derives from the fact that the second expectation on the right-hand-side of \eqref{eq_kl_1} can be expressed as
$$-\mathbb{E}_{f_1}[\log f_2(x;\theta)] \approx -n^{-1}\sum_{i=1}^n\log f_2(x_i;\theta) = -n^{-1}\log L(\theta|x_1,\ldots,x_n),$$
where $L(\theta|x_1,\ldots,x_n)$ is the likelihood function. Therefore, the negative of the average of the log-likelihood is a proxy to obtain the Kullback--Leibler divergence. As such, the minimization of the divergence in \eqref{eq_kl_1} can be obtained, as asymptotic approximation, by maximizing the (expected) likelihood function. That is
$$\min_{\theta}D_{KL}(f_1(x;\eta)\|f_2(x;\theta)) \approx \mathbb{E}_{f_1}[\log f_1(x;\eta)] - \mathbb{E}_{f_1}[\log f_2(x;\hat{\theta)}],$$
with $\hat{\theta}$ being the maximum likelihood estimator of the unknown parameter $\theta$, under the assumption that the true model is $f_1$. Throughout the paper we will assume that the maximum likelihood estimator exists, although not always analytically obtainable.

Let us now assume that both model $f_1$ and $f_2$ are location-scale densities. That is
$$f_1(x;a_1,b_1) = b_1^{-1}h_1\left(\frac{x-a_1}{b_1}\right)\qquad \mbox{and}\qquad f_2(x;a_2,b_2) = b_2^{-1}h_2\left(\frac{x-a_2}{b_2}\right).$$ As shown above, minimizing the Kullback--Leibler divergence in any direction, that is measured either from $f_1$ to $f_2$ or from $f_2$ to $f_1$, corresponds to maximize the likelihood function with respect to the parameters of the second model - i.e. $(a_2,b_2)$ and $(a_1,b_1)$, respectively. In the first case we assume that the data has been generated by $f_1$, whilst in the second case it is $f_2$ the true data-generating model. In \cite{Antle:1969} and, for example, \cite{Dumon:1973} it shown that, for location-scale distributions, the maximum likelihood does not depend from the parameters of the true model. Therefore, the minimum divergence $D_{KL}(f_1(x;a_1,b_1)\|f_2(x;a_2,b_2))$ does not depend on $(a_1,b_1)$ and, similarly, $D_{KL}(f_2(x;a_2,b_2)\|f_1(x;a_1,b_1))$ does not depend on $(a_2,b_2)$.

Furthermore, the property is extended to distributions that can be transformed into location-scale models. In fact, as shown in \cite{DumonAnt:1973}, the maximum likelihood independence from the parameters of the true model is proved to hold for distributions transformable into location-scale models. The result is also derivable by considering the invariance property of the Kullback--Leibler divergence with respect to one-to-one transformations of the variable.

To illustrate the property, we discuss two examples. The first one involves two location-scale distributions: the half-normal and the exponential. The second example involves the Weibull distribution and the lognormal distribution, which can be both transformed into location-scale densities by considering $z = \log x$.

\begin{example}\label{exe1}
Let us consider the half-normal distribution
$$f_1(x;\sigma) = \frac{1}{\sigma}\sqrt{\frac{2}{\pi}}\exp\left(-\frac{x^2}{2\sigma^2}\right) \qquad x\geq0,\sigma>0,$$
and the exponential distribution with the following parametrization
$$f_2(x;\beta) = \frac{1}{\beta}\exp\left(-\frac{x}{\beta}\right)\qquad x\geq0,\beta>0.$$
The Kullback--Leibler divergence from $f_1(x;\sigma^2)$ to $f_2(x;\beta)$, $f_1$ and $f_2$ for simplicity in the notation, is given by
\begin{eqnarray}\label{eq_ex1_1}
D_{KL}(f_1\|f_2) &=& \int f_1\left\{-\log\sigma + \frac{1}{2} \log\left(\frac{2}{\pi}\right) - \frac{x^2}{2\sigma^2}\right\} \; dx \nonumber\\
&& - \int f_1 \left\{-\log\beta-\frac{x}{\beta}\right\}\;dx \nonumber \\
&=& -\log\sigma + \frac{1}{2}\log\left(\frac{2}{\pi}\right)-\frac{\mathbb{E}(x^2)}{2\sigma^2}+\log\beta + \frac{\mathbb{E}(x)}{\beta},
\end{eqnarray}
where the expectations are taken with respect to $f_1$: $\mathbb{E}(x) = \sigma\sqrt{2/\pi}$ and $\mathbb{E}(x^2) = \sigma^2$. The minimum of the divergence \eqref{eq_ex1_1} can be found analytically, and is attained when the two densities have the same mean. Therefore, we set $\beta=\sigma\sqrt{2/\pi}$ to find
$$\min_\beta D_{KL}(f_1\|f_2) = \log(2/\pi)+0.5 = 0.0484.$$
We see that the minimum Kullback--Leibler divergence between $f_1$ and $f_2$ does not depend on $\sigma$. It is possible to see how the result reconciles with the relationship between the minimum divergence and the maximum of the likelihood function. In fact, the maximum likelihood estimator for $\beta$ is the sample mean, which asymptotically tends to the mean, for $f_1$ in this case.
With a similar procedure we find the minimum Kullback--Leibler divergence between $f_2$ and $f_1$. First we have
\begin{eqnarray}\label{eq_ex1_2}
D_{KL}(f_2\|f_1) &=& \int f_2\left\{-\log\beta - \frac{x}{\beta}\right\}\;dx - \int f_2\left\{-\log\sigma+\frac{1}{2}\log\left(\frac{2}{\pi}\right)\right. \nonumber \\
&& \left.-\frac{x^2}{2\sigma^2}\right\}\;dx \nonumber \\
&=& -\log\beta -\frac{\mathbb{E}(x)}{\beta} + \log\sigma - \frac{1}{2}\log\left(\frac{2}{\pi}\right) + \frac{\mathbb{E}(x^2)}{2\sigma^2},
\end{eqnarray}
where, in this case, the expectations are with respect to $f_2$, with values $\mathbb{E}(x)=\beta$ and $\mathbb{E}(x^2) = 2\beta^2$. The minimum of \eqref{eq_ex1_2} is attained for $\sigma=\beta$ and it has value
$$\min_\sigma D_{KL}(f_2\|f_1) = -\frac{1}{2}\log\frac{2}{\pi} = 0.2258$$
As it can be analytically verified, the minimum is obtained when the densities have equal scale parameter. The minimum Kullback--Leibler is independent from $\beta$, as expected. We can again infer the minimization condition from the relationship with the maximum likelihood. Knowing that the maximum likelihood estimator of the parameter $\sigma$ is the sample standard deviation, we have the minimum attained when $\sigma$ equals the standard deviation under $f_2$.
\end{example}

In the next example we compare a lognormal density with a Weibull density, which can be transformed in location-scale parameter models by applying $z=\log x$. In particular, a lognormal distribution with log-scale parameter $\mu$ and shape parameter $c$ can be transformed, by taking the logarithm as above, in a normal with mean $\mu$ and standard deviation $c$. Note that in the example, for convenience, we have replaced the shape parameter by the square of its inverse: $\tau=1/c^2$. Similarly, by taking the logarithm of a Weibull random variable with scale parameter $\lambda$ and shape parameter $\kappa$, we obtain a Gumbel density with location parameter $a=\log\lambda$ and scale parameter $b=1/\kappa$.

\begin{example}\label{exe2}
Let us consider the lognormal density
$$f_1(x;\mu,\tau) = \frac{1}{x}\left(\frac{\tau}{2\pi}\right)^{1/2}\exp\left\{-\frac{\tau}{2}(\log x-\mu)^2\right\}\qquad x>0,\mu>0,\tau>0,$$
and the Weibull density
$$f_2(x;\lambda,\kappa) = \frac{\kappa}{\lambda}\left(\frac{x}{\lambda}\right)^{\kappa-1}\exp\left\{-\left(\frac{x}{\lambda}\right)^\kappa\right\}\qquad x>0,\lambda>0,\kappa>0.$$
The Kullback--Leibler divergence between $f_1(x;\mu,\tau)$ and $f_2(x;\lambda,\kappa)$, $f_1$ and $f_2$ for simplicity in the notation, is given by
\begin{eqnarray}\label{eq_ex2_1}
D_{KL}(f_1\|f_2) &=& \int f_1\left\{-\log(x)-\frac{1}{2}\log\left(\frac{\tau}{2\pi}\right)-\frac{\tau(\log x-\mu)^2}{2}\right\}\;dx \nonumber \\
&&- \int f_1\left\{\log\left(\frac{\kappa}{\lambda}\right)-(\kappa-1)\log\left(\frac{x}{\lambda}\right)-\left(\frac{x}{\lambda}\right)^2\right\}\;dx \nonumber \\
&=& -\mathbb{E}(\log x) - \frac{1}{2}\log\left(\frac{\tau}{2\pi}\right)-\frac{\tau}{2}\mathbb{E}[(\log x - \mu)^2] - \log\left(\frac{\kappa}{\lambda}\right) \nonumber \\
&& - (\kappa-1)\mathbb{E}[\log(x/\lambda)]-\mathbb{E}[(x/\lambda)^2],
\end{eqnarray}
where the expectations are taken with respect to $f_1$. Therefore, we have $\mathbb{E}(\log x) = \mu$, $\mathbb{E}(x^\kappa)=\exp\{k^2/(2\tau)+\mu\kappa\}$ and $\mathbb{E}(\log^2x)=1/\tau+\mu^2$. The minimum divergence in \eqref{eq_ex2_1} is attained at $\lambda=\exp\{1/(2\sqrt{\tau})+\mu\}$ and $\kappa=\sqrt{\tau}$, giving
$$\min_{\lambda,\kappa}D_{KL}(f_1\|f_2) = 1-\frac{1}{2}\log (2\pi) = 0.0811.$$
The minimum Kullback--Leibler between a lognormal and a Weibull distribution, where the minimization is intended with respect to the parameters of the Weibull, does not depend on the parameters of the lognormal $\mu$ and $\tau$. As the maximum of the likelihood function for a Gumbel, with respect to $a=\log\lambda$ and $b=1/\kappa$, cannot be found analytically, it is not possible to show the relationship with the minimum Kullback--Leibler divergence.
The divergence between the Weibull and the lognormal is given by
\begin{eqnarray}\label{eq_ex2_2}
D_{KL}(f_2\|f_1) &=&\int f_2\left\{\log\left(\frac{\kappa}{\lambda}\right)-(\kappa-1)\log\left(\frac{x}{\lambda}\right)-\left(\frac{x}{\lambda}\right)^2\right\}\;dx \nonumber \\
&&- \int f_2\left\{-\log(x)-\frac{1}{2}\log\left(\frac{\tau}{2\pi}\right)-\frac{\tau(\log x-\mu)^2}{2}\right\}\;dx \nonumber \\
&=& \log\left(\frac{\kappa}{\lambda}\right) - (\kappa-1)\mathbb{E}[\log(x/\lambda)] - \mathbb{E}[(x/\lambda)^2] + \mathbb{E}(\log x) \nonumber \\
&&+ \frac{1}{2}\log\left(\frac{\tau}{2\pi}\right) -\frac{\tau}{2}\mathbb{E}[(\log x-\mu)^2].
\end{eqnarray}
The expectations are whit respect to $f_2$: $\mathbb{E}(\log x) = \log\lambda-\gamma/\kappa$, $\mathbb{E}(x^\kappa)=\lambda^\kappa$ and $\mathbb{E}[(\log x-\mu)^2]=6\kappa^2/\pi^2$, where $\gamma$ is the Eulero--Mascheroni constant. The minimum of \eqref{eq_ex2_2} is attained when the two distributions have the same mean and the same variance, that is for $\mu=\mathbb{E}(\log x)$ and $\tau=6(\kappa/\pi)^2$. Therefore
$$\min_{\mu,\tau}D_{KL}(f_2\|f_1) = \frac{1}{2}\log(2\pi)+\log\pi-\gamma-\frac{1}{2}\log6-\frac{1}{2} = 0.0906.$$
The minimum divergence does not depend on the parameters of the Weibull, as expected.  The relationship with the maximum likelihood can be seen if we consider the log-transformations of the distributions, as seen above. In fact, the maximum likelihood estimators of $\mu$ and $\tau$ are, respectively, the sample mean and the reciprocal of the sample variance; asymptotically they become mean and precision of a normal density.
\end{example}

As expected, in general, it is possible to obtain the minimum Kullback--Leibler divergence from the maximum likelihood estimators of the parameter of the model to where the distance is considered. However, as seen in Example \ref{exe2}, the relationship between the minimum Kullback--Leibler divergence and the maximum likelihood appears not to be analytically verifiable when the likelihood function has to be maximized through numerical methods.

\section{Objective model priors for location-scale models}\label{sc_modsel}
In this section, we show how the property discussed in Section \ref{sc_minkl} can be useful in performing Bayesian model selection when the competing models belong to the location-scale family.

Bayes factors represent a well established mean to perform Bayesian model selection, as discussed in \cite{Kass:1995} and \cite{BerPer:2001}, for example. Let us consider models $M_1=\{f_1(x;\theta_1),\pi_1(\theta_1)\}$ and $M_2=\{f_2(x;\theta_2),\pi_2(\theta_2)\}$, where $\pi_1(\theta_1)$ and $\pi_2(\theta_2)$ are the prior distributions for each model-specific parameter (or vector of parameters). If $P(M_i|D)$ is the posterior probability of model $M_i$, for $i=1,2$, given observations $D$, we can define the posterior odds by
\begin{equation}\label{eq_modsel_1}
\frac{P(M_1|D)}{P(M_2|D)} = \frac{P(D|M_1)}{P(D|M_2)}\frac{P(M_1)}{P(M_2)},
\end{equation}
where
$$B_{12} = \frac{P(D|M_1)}{P(D|M_2)},$$
is the Bayes factor. Therefore, the posterior odds are the result of the multiplication of the Bayes factor by the prior odds ($P(M_1)/P(M_2)$), and will give an indication on which model prefer. Note that models $M_1$ and $M_2$ can be replaced by two hypothesis, say $H_1$ and $H_2$, to be tested. In the realm of objective Bayes, when the two models are non-nested, there are not many choices if defining model priors besides the trivial uniform $P(M_1)=P(M_2)=1/2$. However, in \cite{Villa:Walker:2014}, it is proposed a method that takes into consideration the loss in information derived from the choice of the wrong model. The idea is that, if we choose model $M_1$ when $M_2$ is the true model, according to a well known asymptotic Bayesian property \citep{Berk:1966}, the posterior accumulates at the model which is the nearest, in terms of the Kullback--Leibler divergence, to the true model: therefore, $\displaystyle\min_{\theta_2}D_{KL}(f_1(x;\theta_1)\|f_2(x;\theta_2))$, assuming $\theta_1$ known, represents the loss in information in choosing the ``wrong'' model. However, since we do not know $\theta_1$, but we have the prior $\pi_1(\theta_1)$, we can compute the expected loss as
\begin{equation}\label{eq_modsel_2}
\int_{\Theta_1}\min_{\theta_2}D_{KL}(f_1(x;\theta_1)\|f_2(x;\theta_2))\pi_1(\theta_1)\;d\theta_1.
\end{equation}
The model prior is then determined by means of the \emph{self-information} loss function \citep{Merhav:Feder:1998}, which represents the loss connected to a probability statement. For example, the self-information loss for model $M$ is $-\log P(M)$. Therefore, by equating the self-information loss with the expected loss in \eqref{eq_modsel_2}, we have that the prior for $M_1$ is
\begin{equation}\label{eq_modsel_3}
P(M_1) \propto \exp\left\{\int_{\Theta_1}\min_{\theta_2}D_{KL}(f_1(x;\theta_1)\|f_2(x;\theta_2))\pi_1(\theta_1)\;d\theta_1\right\},
\end{equation}
and the prior for $M_2$ is
\begin{equation}\label{eq_modsel_4}
P(M_2) \propto \exp\left\{\int_{\Theta_2}\min_{\theta_1}D_{KL}(f_2(x;\theta_2)\|f_1(x;\theta_1))\pi_2(\theta_2)\;d\theta_2\right\}.
\end{equation}

When the models $M_1$ and $M_2$ are location-scale distributions (or transformable into members of the family), we have shown in Section \ref{sc_minkl} that the minimum Kullback--Leibler divergence, with respect to the parameters of the second model, does not depend on the parameters of the first model. Given that the model priors \eqref{eq_modsel_3} and \eqref{eq_modsel_4} are based on minimum Kullback--Leibler divergences, we see that, if $M_1$ and $M_2$ are location-scale models, the prior odds in \eqref{eq_modsel_1} will not depend on the choice of priors $\pi_1(\theta_1)$ and $\pi_2(\theta_2)$.  From an objective point of view, this result means that the prior assigned to each model depend only on the information conveyed by the choice of models; and it represents a less informative approach than the one of assigning equal model probability \textit{a priori}. \\

Let us see how the above result can be implemented by considering the two examples introduced in Section \ref{sc_minkl}.

In Example \ref{exe1} we compare a half-normal distribution with an exponential distribution. To compute the posterior odds, in order to assess which model better represents the observed data, we consider the prior odds given by the ratio between $P(M_1)\propto1.05$ and $P(M_2)\propto1.25$ (or vice versa). Normalizing, we have $P(M_1)=0.46$ and $P(M_2)=0.54$. We see that the loss in information in considering the \emph{wrong} model is not necessarily symmetrical. That is, for the specific example, if we choose the exponential model when the true one is the half-normal, we have (\textit{a priori}) a larger loss than in the opposite scenario. In addition, there is a clear difference in assigning model prior accordingly to this approach compared to the uniform prior.

However, the result is not common to all the cases. In fact, if we consider Example \ref{exe2}, we can see that $P(M_1)\propto1.08$ and $P(M_2)\propto1.09$.  The normalized model prior probabilities, $P(M_1)=P(M_2)=0.50$, are uniformly distributed.

\section{Discussion}\label{sc_disc}
We consider the loss in information deriving from the choice of a model that approximates the unknown true model, where the models belong to the family of location-scale distributions. In statistical inference the aim is to have an approximated model that is as close as possible to the true one; in other words, the most similar. As the dissimilarity between to densities is naturally measured by the Kullback--Leibler divergence, it goes that the aim is to choose the approximate model that minimizes the divergence from the true model.

In this paper, we show that the minimum Kullback--Leibler divergence between two location-scale models does not depend on the parameters of the model from where the divergence is measured. In addition, the property holds for distributions that can be transformed into location-scale models.

A possible implication of the property is in objective Bayesian model selection, where a prior on the model space can be defined by taking into consideration the expected loss in information that the choice of the ``wrong'' model produces. Interestingly, we show that the loss in information is not necessarily symmetrical, and that some location-scale models have \emph{a priori} a relatively higher importance compared to others. If we consider the subset of genuine location-scale distributions, we have the following result. Both the model prior and the parameter priors can always be objective. As objective priors are in general improper, their use in model selection by means of Bayes factors is limited to parameters that are common to both models; therefore, applicable to location-scale model selection problems. By assigning model priors on the basis of the information carried by the choice of the competing models, the Bayesian procedure can be seen as objective as it can get, given that a uniform prior on the model space can arguably be considered as noninformative.

Surely, there may be other scenarios where the property of the minimized Kullback--Leibler divergence for location-scale models here discussed could be useful. However, in this paper we have limited our considerations to objective Bayesian model selection.


\begin{thebibliography}{}
\bibitem[\protect\astroncite{Akaike}{1973}]{Aka:1973}
{\sc Akaike, H.} (1973).
\newblock Information theory as an extension of the maximum likelihood principle.
\newblock In {\em Second International Symposium on Information Theory} (eds. B.~N. Petrov, and F. Csaki), Akademiai Kiado, Budapest, 267--281.

\bibitem[\protect\astroncite{Antle and Bain}{1969}]{Antle:1969}
{\sc Antle, C.~E. and Bain, L.~J.} (1969).
\newblock A property of maximum likelihood estimators of location and scale parameters.
\newblock {\em SIAM Review} {\bf 11}, 251--253

\bibitem[\protect\astroncite{Berger and Pericchi}{2001}]{BerPer:2001}
{\sc Berger, J.O. and Pericchi, L.~R.} (2001).
\newblock Objective Bayesian methods for model selection: introduction and comparison.
\newblock {\em IMS Lecture Notes -- Monograph Series} {\bf 38}, 135--193

\bibitem[\protect\citeauthoryear{Berk}{1966}]{Berk:1966}
{\sc Berk, R.~H.} (1966).
\newblock Limiting behaviour of posterior distributions when the model is incorrect.
\newblock {\em Ann. of Math. Statist.} {\bf 37}, 51--58.

\bibitem[\protect\astroncite{Dumonceaux et al.}{1973a}]{Dumon:1973}
{\sc Dumonceaux, R., Antle, C.~E. and Haas, G.} (1973a).
\newblock Likelihood ratio test for discrimination between two models with unknown location and scale parameters.
\newblock {\em Technometrics} {\bf 15}, 19--27

\bibitem[\protect\astroncite{Dumonceaux and Antle}{1973b}]{DumonAnt:1973}
{\sc Dumonceaux, R. and Antle, C.~E.} (1973b).
\newblock Discrimination between the log-normal and the Weibull distributions.
\newblock {\em Technometrics} {\bf 15}, 923--926

\bibitem[\protect\astroncite{Johnson et al.}{1994}]{John:1994}
{\sc Johnson, N.~L., Kotz, S. and Balakrishan, N.} (1994).
\newblock {\em Continuous Univariate Distributions}.
\newblock Vol. 1 and 2, 2nd ed., Wiley, New York

\bibitem[\protect\astroncite{Johnson et al.}{1992}]{John:1992}
{\sc Johnson, N.~L., Kotz, S. and Kemp, A.} (1992).
\newblock {\em Univariate Discrete Distributions}.
\newblock 2nd ed., Wiley, New York

\bibitem[\protect\astroncite{Kass and Raftery}{1995}]{Kass:1995}
{\sc Kass, R.~R. and Raftery, A.~E.} (1995).
\newblock Bayes factors.
\newblock {\em J. Amer. Statist. Assoc.} {\bf 90}, 773--795

\bibitem[\protect\astroncite{Kullback and Leibler}{1951}]{Kull:1951}
{\sc Kullback, S and Leibler, R.~A.} (1951).
\newblock On information and sufficiency.
\newblock {\em Annals of Mathematical Statistics} {\bf 22}, 79--86.

\bibitem[\protect\astroncite{Lv and Liu}{2014}]{LvLiu:2014}
{\sc Lv, J. and Liu, J.~S.} (2014).
\newblock Model selection principles in misspecified models.
\newblock {\em J. R. Statist. Soc. B} {\bf 76}, 141--167

\bibitem[\protect\citeauthoryear{Merhav and Feder}{1998}]{Merhav:Feder:1998}
{\sc Merhav, N. and Feder, M.} (1998).
\newblock Universal prediction.
\newblock {\em IEEE Trans. Inf. Theory} {\bf 44}, 2124--2147.

\bibitem[\protect\citeauthoryear{Villa and Walker}{2014}]{Villa:Walker:2014}
{\sc Villa C. and Walker, S.~G.} (2014).
\newblock An objective Bayesian criterion to determine model prior probabilities.
\newblock {\em Scandinavian Journal of Statistics} {\bf 42}, 947--966

\bibitem[\protect\astroncite{Zheng et al.}{2004}]{Zheng:2004}
{\sc Zheng, G, Freidlin, B., and Gaswirth, J.~L.} (2004).
\newblock Using Kullback--Leibler information for model selection when the data-generating model is unknown: applications to genetic testing problems.
\newblock {\em Technometrics} {\bf 15}, 19--27
\end{thebibliography}


\end{document}